\documentclass[11pt]{amsart}

\usepackage{amssymb}

\evensidemargin\oddsidemargin

\begin{document}
\baselineskip=18pt
\setcounter{page}{1}
    
\newtheorem{Conju}{Conjecture 1\!\!}
\newtheorem{Conjd}{Conjecture 2\!\!}
\newtheorem{Theo}{Theorem\!\!}
\newtheorem{Rqs}{Remarks}
\newtheorem{Lemm}{Lemma\!\!}
\newtheorem{Def}{Definition\!\!}
\newtheorem{Prop}{Proposition\!\!}

\renewcommand{\theConju}{}
\renewcommand{\theConjd}{}
\renewcommand{\theDef}{}
\renewcommand{\theTheo}{}
\renewcommand{\theLemm}{}
\renewcommand{\theProp}{}

\def\a{\alpha}
\def\b{\beta}
\def\B{{\bf B}} 
\def\C{{\mathcal{C}}} 
\def\CC{{\mathbb{C}}} 
\def\E{{\mathcal{E}}} 
\def\Da{{\rm D}_\a}
\def\Dq{{\rm D}_q}
\def\EE{{\mathbb{E}}} 
\def\elaw{\stackrel{d}{=}}
\def\eps{\varepsilon}
\def\F{{\bf F}} 
\def\G{\gamma} 
\def\fa{f_\a} 
\def\HH{{\mathbb{H}}} 
\def\hS{{\hat S}}
\def\hT{{\hat T}}
\def\hX{{\hat X}}
\def\ii{{\rm i}}
\def\K{{\bf K}} 
\def\L{{\mathcal{L}}} 
\def\lb{\lambda}
\def\lacc{\left\{}
\def\lcr{\left[}
\def\lpa{\left(}
\def\lva{\left|}
\def\NN{{\mathbb{N}}} 
\def\pb{{\mathbb{P}}}
\def\R{{\mathcal{R}}}
\def\rl{{\mathbb{R}}}
\def\racc{\right\}}
\def\rcr{\right]}
\def\rpa{\right)}
\def\rva{\right|}
\def\T{{\bf T}} 
\def\TT{{\rm T}} 
\def\U{{\bf U}} 
\def\Un{{\bf 1}}
\def\Xa{X_\a} 
\def\Ya{Y_\a} 
\def\Yan{Y_{\a,n}} 
\def\ZZ{{\mathbb{Z}}} 

\newcommand{\fin}{\vspace{-0.4cm}
                  \begin{flushright}
                  \mbox{$\Box$}
                  \end{flushright}
                  \noindent}

\title[Positive stable densities and bell-shape]{Positive stable densities and the bell-shape}

\author[Thomas Simon]{Thomas Simon}

\address{Laboratoire Paul Painlev\'e, Universit\'e Lille 1, Cit\'e Scientifique, F-59655 Villeneuve d'Ascq Cedex. {\em Email}: {\tt simon@math.univ-lille1.fr}}



\begin{abstract} We show that positive stable densities are bell-shaped, that is their $n$-th derivatives vanish exactly $n$ times on $(0,+\infty)$ and have an alternating sign sequence. This confirms the graphic predictions of Holt and Crow (1973) in the positive case.
\end{abstract}

\maketitle

\section{Introduction}

Consider a smooth probability density defined on some open interval $I \subset\rl$ and such that all its derivatives vanish at both ends of $I$.
Rolle's theorem shows that the $n-$th derivative vanishes at least $n$ times on $I.$ The density is said to be bell-shaped if the $n-$th derivative vanishes exactly $n$ times on $I$ for every $n\ge 1.$ For $n =1$ this amounts to strict unimodality. For $n=2$ this means, as for the familiar bell curve, that there is one inflection point on each side of the mode and that the second derivative is successively positive, negative, and positive. Rolle's theorem also entails that the zeroes of the successive derivatives of a bell-shaped density strictly interlace, and hence that each derivative has an alternating sign sequence, starting positive. The sign sequence of the third derivative can be observed on a curve in comparing the aberrancy and the slope, as shown by Transon's formula \cite{Sc}, but contrary to the second derivative it requires an educated eye to guess what this sequence must be at first glance. The fourth derivative of a density has also a kinematic interpretation in terms of penosculating conics - see \cite{Sc} and the references therein for a complete account, which however cannot be understood through the sole sequence of its signs. 

For a given density whose derivatives vanish at both ends of its definition interval, the usual way to show the bell-shape is to factorize its $n$-th derivative by a positive function and a polynomial function of degree $n$. This entails that this $n$-th derivative vanishes at most, and hence exactly, $n$ times. Basic examples - whose details can be checked by the reader and serve as an exercise in calculus for sophomores - are the Gaussian, Gumbel and centered Student densities on $\rl$, and the inverse Gamma densities on $(0, +\infty).$ When the density is not explicit, such factorizations are not appropriate anymore. Another criterion for the bell-shape is the ETP character of the associated additive convolution kernel - see Chapter 6.11.C in \cite{K} - but the latter is rather difficult to check for non-explicit densities, and also quite stringent since it is e.g. fulfilled neither by the Student nor by the inverse Gamma densities. In this paper we are interested in the positive $\a-$stable random variables $\Xa$ and their densities $\fa, \, 0<\a<1.$ The latter are characterized by their Laplace transform and we choose the normalisation
$$\EE[e^{-\lb \Xa}]\; =\;\int_0^\infty e^{-\lb x} \fa(x)\, dx\; =\; e^{-\lb^\a}, \quad \lb \ge 0.$$
It is known - see e.g. the expansions (14.31) and (14.35) in \cite{S} - that $\fa$ is real-analytic on $(0,+\infty)$ and that all its derivatives vanish at zero and at infinity. By Proposition 7.1.3 in \cite{K} and (14.31) in \cite{S}, the kernel $\fa(x-y)$ is not ${\rm TP}_2$ and hence not ETP. The function $f_{1/2}$ is explicit thanks to the identity in law
$$X_{1/2}\;\elaw\; \frac{1}{4\Gamma_{1/2}}$$
where here and throughout $\Gamma_a$ denotes the standard Gamma random variable of parameter $a > 0,$ and is hence bell-shaped. The functions $f_{1/3}$ and $f_{2/3}$ can be written down in terms of a modified Bessel resp. a confluent hypergeometric function - see e.g.  (2.8.31) and (2.8.33) in \cite{Z} - but it does not seem to the author that the associated second order ODE's can provide any substantial information on their bell-shape. 


The strict unimodality of $\fa$ was first proved by Ibragimov and Chernin - see the proof of the Theorem in \cite{IC} - and then extended to all real stable densities by Sato and Yamazato - see Theorem 1.4 in \cite{SY}. Since the latter are real-analytic, this actually follows from Yamazato's original theorem \cite{Y1} on self-decomposable distributions. See also \cite{S2} for a short proof. The graphic simulations made by Holt and Crow \cite{HC}, partly reproduced in \cite{Z}, display bell-shaped curves for all real stable densities, at least at the level of the second derivative. A consequence of Theorem 5.1.(ii) in \cite{SY} is that $\fa$ has only one inflection point on the left side of its mode, and more generally that each $n$-th derivative of $\fa$ vanishes only once before the first zero of the $(n-1)$-th derivative. The bell-shape property was claimed in \cite{G1} for all stable densities but a serious mistake on the TP character of the kernel $\fa(xy^{-1})$ - see Remark (f) in \cite{S2} for an explanation - invalidates this result, also at the unimodal level, except for a particular class of stable densities which will be discussed in Section 3.3 and does not include the positive case. Apart from \cite{HC, S2, G1} the bell-shape of stable densities seems to have escaped investigation, even in the most visual case $n=2.$ In this note, we show the full property in the positive case. 
 
\begin{Theo} The densities $\fa$ are bell-shaped.
\end{Theo}

Our main argument comes from Schoenberg's variation-diminishing property and the total positivity of certain infinitely convoluted exponential kernels. The little known fact that $\Xa$ admits infinite exponential sums as additive factors had been noticed in \cite{Y2}. It is used here in conjunction with a precise analysis of the other factor, which contains a certain exponential mixture with no atom at zero and can be chosen in order to have a weak bell-shape property of arbitrary order. In the case when $1/\a$ is an integer, another proof of the bell-shape is obtained quickly by a multiplicative factorization with inverse Gamma laws - see Section 3.1 below, but the argument does not extend to the other cases. With this method, the bell-shape of the density of $\log X_{1/n}$ on $\rl$ can also be established. In Section 3.2, we raise two natural conjectures on the bell-shape of positive self-decomposable densities.

 \section{Proof of the theorem}

We first establish an additive factorization of $\Xa$ via standard exponential variables. Recall that a positive random variable $X\sim\mu$ is said to be an exponential mixture (ME) if its law has the form
$$\mu(dx)\; =\; c\delta_0(dx)\; +\; f(x) dx$$
with $c\in [0,1]$ and $f$ a completely monotone function on $(0,+\infty).$ When $c=0$ we will use the notation $X\in{\rm ME}^*,$ and in this case Bernstein's theorem shows that the density $f$ writes
\begin{equation}
\label{ME}
f(x) \; =\; \int_0^\infty \theta e^{-\theta x} \mu(d\theta)
\end{equation}
over $(0, +\infty),$ where $\mu$ is some probability measure on $(0,+\infty).$ In other terms $X\elaw {\rm Exp}\, (1)\times X_\mu^{-1}$ 
where here and throughout Exp $(\lb)$ stands for the exponential law with parameter $\lb,$ and $X_\mu \sim\mu.$ The following lemma is mostly due to Yamazato  - see \cite{Y2} pp. 601-602 - but we rephrase it and give a proof for the reader's convenience.

\begin{Lemm}[Yamazato] One has the independent factorization
\begin{eqnarray}
\label{Dixmude}
\Xa\;\elaw\;\Ya\; +\; \sum_{n\ge 1} \,\Yan,
\end{eqnarray}
where $\Yan\sim{\rm Exp} ((n\pi/\sin(\pi\a))^{1/\a})$ for every $n\ge 1$ and $\Ya\in{\rm ME}^*.$ 
\end{Lemm}

\proof We first write
$$\lb^\a\; =\; \frac{\a}{\Gamma(1-\a)}\int_0^\infty (1-e^{-\lb x})\frac{dx}{x^{\a+1}}\; =\; \int_0^\infty (1-e^{-\lb x})\lpa\int_0^\infty (c_\a u)^\a e^{-xu} du\rpa dx$$
with the notation $c_\a = \sin(\pi\a)/\pi.$ This yields the further decomposition
$$\lb^\a\; =\;\int_0^\infty (1-e^{-\lb x})l_\a(x)dx\; +\; \int_0^\infty (1-e^{-\lb x})\lpa\int_0^\infty \lcr (c_\a u)^\a\rcr e^{-xu} du\rpa dx$$
where $[.]$ stands for the integer part and
$$l_\a(x)\; =\; \int_0^\infty ( (c_\a u)^\a -\lcr (c_\a u)^\a\rcr) e^{-xu} du$$ 
is a completely monotone function. By Theorem 51.10 in \cite{S} there exists a random variable $\Ya$ belonging to the Bondesson class such that
$$-\log\EE[e^{-\lb \Ya}]\; =\; \int_0^\infty (1-e^{-\lb x})l_\a(x)dx, \qquad \lb \ge 0.$$ 
Because $0\le (c_\a u)^\a -\lcr (c_\a u)^\a\rcr\le 1$ for all $u\ge 0$ and $u^{-1}((c_\a u)^\a -\lcr (c_\a u)^\a\rcr)$ is integrable at zero, Steutel's criterion - see Theorem 51.12 in \cite{S} - shows that $\Ya$ is actually ME. We next prove that $\Ya$ has no atom at zero, which will entail $\Ya\in{\rm ME}^*.$ A change of variable and an integration by parts yield
$$l_\a(x)\; =\; \frac{x}{c_\a}\int_0^\infty e^{-xc_\a^{-1}u}  \lpa\int_0^u (t^\a -\lcr t^\a\rcr)\, dt\rpa du$$
for every $x>0,$ and it is a bit tedious but elementary to see that
$$\int_0^u (t^\a -\lcr t^\a\rcr)\, dt\; \sim\; \frac{[u^\a]^{1/\a}}{2}\; \sim\; \frac{u}{2}$$
as $u\to+\infty.$ Theorem XIII.5.4 in \cite{F} entails
\begin{equation}
\label{Tauber}
l_\a(x)\; \sim\;\frac{c_\a}{2x}
\end{equation}
as $x\to 0,$ so that $\Ya$ has infinite L\'evy measure and, by Theorem 27.4 in \cite{S}, no atom at zero. On the other hand, setting $\kappa_{\a,n} = (n\pi/\sin(\pi\a))^{1/\a}$ for every $n\ge 1,$ one has
$$\int_0^\infty \lcr (c_\a u)^\a\rcr e^{-xu} du\; =\; \sum_{n\ge 1} \int_{\kappa_{\a,n}}^\infty e^{-xu} du\; =\; \sum_{n\ge 1} \frac{e^{-x\kappa_{\a,n}}}{x}\cdot$$
Hence,
\begin{eqnarray*}\int_0^\infty (1-e^{-\lb x})\lpa\int_0^\infty \lcr (c_\a u)^\a\rcr e^{-xu} du\rpa dx & = & \sum_{n\ge 1}\int_0^\infty (1-e^{-\lb x}) \frac{e^{-x\kappa_{\a,n}}}{x}dx\\
& = & \sum_{n\ge 1} (\log(\lb+ \kappa_{\a,n}) - \log(\kappa_{\a,n}))\\
& = & -\log \EE[ \exp [ -\lb\sum_{n\ge 1} \Yan]],
\end{eqnarray*}
where the first equality comes from Fubini's theorem and the second from Frullani's. Putting everything together completes the proof.
 
\endproof

\begin{Rqs} {\em (a) Because $\a\in(0,1),$ the sum on the right-hand side of (\ref{Dixmude}) has finite expectation. This shows that $\Xa$ and $\Ya$ have the same $\L_p-$integrability index which is $\a,$ excluded. The latter can also be seen from Theorem 25.3 in \cite{S} and the fact that $l_\a(x)\sim (\a/\Gamma(1-\a)) x^{-(\a +1)}$ as $x\to +\infty.$ 

(b) Since the infinite exponential sum has a ${\rm PF}_2$ viz. log-concave density - see the proof of the theorem below for a stronger property - and since the density of $\fa$ is not ${\rm PF}_2,$ one can interpret $\Ya$ as the factor which breaks down the ${\rm PF}_2$ property for $\fa$. Recall indeed that $\fa$ is log-concave until its second inflection point - see Theorem 1.3 (vii) in \cite{SY}. On the other hand $\Ya$ has a log-convex density by H\"older's inequality. It is plausible that $\log\fa$ has only {\em one} inflection point, which can be readily checked for $\a =1/2$.}

\end{Rqs}

We now introduce our notations concerning the sequence of signs of a real function defined on $(0,+\infty).$ We set
$${\rm sign}(x)\; =\; \lacc \begin{array}{cl}+ & \mbox{if $x>0$}\\
0 & \mbox{if $x=0$}\\- & \mbox{if $x<0$}\end{array}\right.\qquad\mbox{and}\qquad \pm^n\; =\; \lacc \begin{array}{cl}+ & \mbox{if $n$ is even}\\
- & \mbox{if $n$ is odd}\end{array}\right.$$
for every $x\in[-\infty,+\infty], n\in\NN.$ For $\{u_n\} = \{\eps_1, \ldots, \eps_n\}$ some finite sequence in $\{-,0,+\},$ we say that a continuous function $f :\, (0,+\infty) \to\rl$ is of type $\eps_1\ldots\eps_n$ if it has limits (finite or infinite) at zero and at infinity, vanishes on a finite set, and if the ordered sequence of its signs on $[0,+\infty]$ is given by $\{\eps_1, \ldots, \eps_n\}.$ For brevity, we will write $f\sim\{ u_n\}$ or $f\sim\eps_1\ldots\eps_n$ to express this property. Observe that by the intermediate value theorem, the zero set of $f\sim\eps_1\ldots\eps_n$ corresponds to a subsequence of period 2 in $\{\eps_1, \ldots, \eps_n\}.$ In particular, a function of type $\eps_1 \cdots \eps_n$ vanishes exactly $(n/2 -1)$ times on $(0,+\infty)$ if $n$ is even and either $(n-1)/2$ or $(n-3)/2$ times on $(0,+\infty)$ if $n$ is odd. For example, the $n$-th derivative of a completely monotone function is of type $\pm^n0$ for every $n\ge 0.$ The density function $\fa$ is of type $0\!+\!0$ and, by the strict unimodality of $\Xa,$ its derivative $\fa'$ is of type $0\!+\!0\!-\!0.$ Introducing the sequences
$$\{a_n\}\;=\; 0\!\pm^0\!0\!\pm^1\!0\cdots\pm^n\!0\qquad \mbox{and} \qquad\{b_n\} = \pm^00\!\pm^1\!0\cdots\pm^{n}\!0$$
for every $n\ge 1,$ we can now state the central definition of this paper.

\begin{Def} For every $n\ge 0$, a smooth function $f:(0,+\infty)\to\rl$  is said to be weakly bell-shaped of order $n$ {\em (${\rm WBS}_n$)} if $f^{(i)}\sim \{a_i\}$ $\forall\,i= 0\ldots n$ and $(-1)^{n+1+i}f^{(i)}\sim \{b_{n+1}\}$ $\forall\, i\ge n+1.$
\end{Def}

The bell-shape property of a function $f$ means that $f^{(i)}\sim \{a_i\}$ $\forall\,i\ge 0,$ so that weakly bell-shaped functions are never bell-shaped. Notice also that ${\rm WBS}_p\,\cap\,{\rm WBS}_q=\emptyset$ if $p\neq q.$ The following result provides our key-argument.

\begin{Prop} Let $X\in {\rm ME}^*$ and $0<\lb_1\le\lb_2\le\ldots\le\lb_n\le\ldots$ be some sequence. For every $n\ge 1$, the independent sum 
$X\,+\, {\rm Exp}(\lb_1) \,+\, \cdots \,+\, {\rm Exp}(\lb_n)$
has a ${\rm WBS}_{n-1}$ density.
\end{Prop}

\proof Set $f_0$ for the density of $X$ and recall from (\ref{ME}) that it is completely monotone. In particular $f_0$ is real-analytic on $(0,+\infty).$ We denote by $f_n$ the density of $X+ {\rm Exp}(\lb_1) + \cdots + {\rm Exp}(\lb_n)$ for every $n\ge 1.$ The latter are connected to one another by the formula
\begin{equation}
\label{Convol}
f_n(x)\; =\; \lb_ne^{-\lb_n x}\int_0^x e^{\lb_n y} f_{n-1}(y) dy,
\end{equation}
and hence solutions on $(0,+\infty)$ to the linear ODE's
\begin{equation}
\label{ODE}
f_n'\; +\; \lb_n f_n\; =\; \lb_n f_{n-1}.
\end{equation} 
An induction hinging on (\ref{Convol}) and (\ref{ODE}) show that all densities $f_n$ are real-analytic on $(0, +\infty),$ that all derivatives $f_n^{(i)}$ vanish at infinity, and that $f_n^{(i)}(0+) = 0$ for all $n\ge i+1\ge 1.$ In particular one has $f_n \sim \{a_0\}$ for every $n\ge 1.$ Using the complete monotonicity of $f_0$ for the initial step, another induction shows $(-1)^{n+i}f_n^{(i)}(0+) > 0$ for all $i\ge n\ge 0,$ possibly with infinite values. This entails $f_n^{(n)}(0+) > 0$ for all $n\ge 0,$ so that $0+$ is an isolated zero of $f_n^{(i)}$ for all $n\ge i+1\ge 1$ (notice that this cannot be proved directly by analyticity because $0+$ is also a singular zero). The weak bell-shape property will now be established through a more elaborate induction on $n$.\\

We first show that $f_1$ is ${\rm WBS}_0.$ Since $f_1 \sim \{a_0\},$  we know that $f_1'$ vanishes at least once on $(0,+\infty).$ Set $x_1 =\inf\{x > 0, \, f_1'(x) =0\} > 0,$ the strict inequality coming from $f_1'(0+) > 0.$ Differentiating (\ref{ODE}) yields 
$$f_1''(x) \; =\; \lb_1f_0'(x)\; <\; 0$$
for every $x >0$ such that $f_1'(x) = 0.$ This shows that $f_1'$ vanishes only at $x_1,$ whence $f_1'\sim\{b_1\}.$ Again, this entails that $f_1''$ vanishes at least once on $(0,+\infty).$ Set $x_2 =\inf\{x > 0, \, f_1''(x) =0\} > 0,$ the strict inequality coming from $f_1''(0+) < 0.$ Differentiating (\ref{ODE}) further entails 
$$f_1'''(x) \; =\; \lb_1f_0''(x)\; >\; 0$$
for every $x >0$ such that $f_1''(x) = 0.$ This shows that $f_1''$ vanishes only at $x_2,$ whence $-f_1''\sim\{b_1\}.$ Using the complete monotonicity of $f_0,$ the same argument yields $(-1)^{i+1}f_1^{(i)}\sim\{b_1\}$ for every $i\ge 1,$ as required.\\

We next show the induction step and suppose that $f_n$ is ${\rm WBS}_{n-1}$ for some $n\ge 1.$ We have already seen that $f_{n+1}\sim\{a_0\}$ and we will first prove by induction on $i$ that $f_{n+1}^{(i)}\sim \{a_i\}$ for every $1\le i\le n$. We begin with $f_{n+1}'$ for the sake of clarity. Since $f_{n+1}\sim\{a_0\},$ we know that $f_{n+1}'$ vanishes at least once on $(0,+\infty)$. Set $x_1 =\inf\{ x> 0, \; f_{n+1}'(x) = 0\}> 0,$ where the strict inequality comes from the fact that $0+$ is an isolated zero of $f_{n+1}'.$ Differentiating (\ref{ODE}) as above entails
$$ \lb_n f_n'(x_1) \; =\; f_{n+1}''(x_1)\; \le\; 0,$$
where the inequality comes from the initial profile $0\!+\!0$ of $f_{n+1}'.$ Suppose that $f_{n+1}''(x_1)= 0.$ Then by analyticity $x_1$ is an isolated zero of $f_{n+1}'',$ so that
$$\lim_{x\uparrow x_1}\frac{f_{n+1}''(x)}{f_{n+1}'(x)}\; =\; -\infty.$$
By (\ref{ODE}), this entails $f_n'(x) < 0$ as $x\uparrow x_1.$ But the induction hypothesis yields $f_n'\sim \{b_1\}$ if $n=1$ and $f_n'\sim \{a_1\}$ if $n >1,$ so that one cannot have $f_n'(x_1) = 0,$ a contradiction. Hence $f_{n+1}''(x_1)< 0$ and if there exists $x_2 > x_1$ such that $f_{n+1}'(x_2) = 0,$ then again the induction hypothesis entails $f_n'(x_2) < 0$ so that $f_{n+1}''(x_2)< 0,$ a contradiction. This shows that $f_{n+1}'$ vanishes only once on $(0,+\infty)$ and, all in all, that $f_{n+1}'\sim\{a_1\}.$

We now suppose $f_{n+1}^{(i)}\sim \{a_i\}$ for some $i< n$ and show $f_{n+1}^{(i+1)}\sim \{a_{i+1}\}.$ We know that $f_{n+1}^{(i+1)}$  vanishes at least $i+1$ times on $(0,+\infty),$ and also from the above that $f_{n+1}^{(i+1)}$ has initial profile $0\!+\!0$. Set $\{x_j, \; j\ge 1\}$ for the ordered sequence of its zeroes on $(0,+\infty).$ Again, one gets from (\ref{ODE})
\begin{eqnarray}
\label{Induc}
f_{n+1}^{(i+2)}(x_j)\, =\, \lb_n f_n^{(i+1)}(x_j)\quad\mbox{and}\quad \lim_{x\uparrow x_j}\frac{f_{n+1}^{(i+2)}(x)}{f_{n+1}^{(i+1)}(x)}\, =\,\lim_{x\uparrow x_j}\frac{f_n^{(i+1)}(x)}{f_{n+1}^{(i+1)}(x)}\, =\, -\infty
\end{eqnarray}
for every $j\ge 1.$ For clarity we consider two separate cases. 

(i) If $i+1=n,$ the induction hypothesis entails $f_n^{(i+1)}\sim\{b_n\} =\pm^0 0\cdots\pm^n\!0.$ Set $I^k_n\leftrightarrow\{\pm^k0\}$ for the $k$-th half-closed interval corresponding to each part $\pm^k0$ of the graph of $f_n^{(i+1)}, k =0\ldots n.$ The last equality in (\ref{Induc}) with $j=1$ and the initial profile $0\!+\!0$ of $f_{n+1}^{(i+1)}$ entail that $x_1\in I^{k_1}_n$ with $k_1$ odd. If $f_n^{(i+1)}(x_1) < 0,$ then (\ref{Induc}) with $j=1,2$ and the induction hypothesis show that $x_2\in I^{k_2}_n$ with $k_2$ even. If $f_n^{(i+1)}(x_1) = 0,$ then $x_2\in I^{k_2}_n$ with $k_2 > k_1.$ Iterating the procedure shows that there are at most, and hence exactly, $n=i+1$ zeroes and that $x_k\in I^k_n$ for every $k=1\ldots n.$ All in all this shows that $f_{n+1}^{(i+1)}\sim\{a_{i+1}\},$ as required. 

(ii) If $i+1<n,$ one has $f_n^{(i+1)}\sim\{a_{i+1}\} =0\!\pm^0\!0\cdots\pm^{i+1}\!0$ by the induction hypothesis. But adding a $0$ at $0+$ in the profile does not change anything in the above analysis, and the fact that $f_{n+1}^{(i+1)}\sim\{a_{i+1}\}$ follows exactly along the same lines as in (i).\\

We finally show that $(-1)^{n+1+i}f_{n+1}^{(i)}\sim \{b_{n+1}\}$ for every $i\ge n+1.$ We only sketch the proof since the arguments are the same. Consider first the case $i=n+1.$ Since $f_{n+1}^{(n)}\sim \{a_n\},$ we know that $f_{n+1}^{(n+1)}$ vanishes at least $n+1$ times on $(0,+\infty)$ and set $\{x_k, \, k\ge 1\}$ for the ordered sequence of its zeroes. Recall $f_{n+1}^{(n+1)}(0+) > 0$ and, from the induction hypothesis, that $f_{n}^{(n+1)}\sim\pm^10\cdots\pm^{n+1}\!0.$ Set again $I^k_n\leftrightarrow\{\pm^k0\}$ for the $k$-th half-closed interval corresponding to each part $\pm^k0$ of the graph of $f_n^{(n+1)}, k =1\ldots n+1.$ The same reasoning as above shows that necessarily $x_k\in I^k_n$ for all $k$ and hence that $f_{n+1}^{(n+1)}$ vanishes exactly $n+1$ times on $(0,+\infty),$ in other terms $f_{n+1}^{(n+1)}\sim \{b_{n+1}\}.$ It is then easy to show by induction on $i$ that $(-1)^{n+1+i}f_{n+1}^{(i)}\sim \{b_{n+1}\}$ for every $i\ge n+1,$ in using $(-1)^{n+1+i}f_{n+1}^{(i)}(0+) >0$ and the induction hypothesis $(-1)^{n+1+i}f_{n}^{(i)}\sim \pm^10\cdots\pm^{n+1}\!0.$

\endproof

\begin{Rqs}{\em (a) The complete monotonicity of $f_0$ is crucial in the above argument, initializing all induction steps. It seems difficult to extend the proposition to broader classes than ME.

(b) When $X$ itself has an exponential law and all parameters are different, the density of $X\,+\, {\rm Exp}(\lb_1) \,+\, \cdots \,+\, {\rm Exp}(\lb_n)$ is a linear combination of $n+1$ distinct exponentials functions. Since the latter define a Chebyshev system - see \cite{K} p. 24, the property ${\rm WBS}_n$ is much easier to prove. We leave the details to the reader.

(c) The method of examining the sign of the $(n+1)$-th derivative at the zero of the $n$-th derivative is used in \cite{IC} for $n=1$ in order to show the strict unimodality of $\fa,$ after some explicit contour integration. The involved computations to extend this argument directly on the further derivatives of $\fa$ seem however extremely heavy, and the author does not believe that they could yield the full bell-shape. See \cite{G2} for an application of the contour integration of \cite{IC} to the asymptotic behaviour of these further derivatives.}
\end{Rqs}

\noindent
{\bf End of the proof.} Fix $n\ge 1$ and set $f_{\a,n}$ resp. $g_{\a,n}$ for the density of
$$\Xa\, +\, \sum_{k=1}^{n+2} \,Y_{\a,k}\qquad\mbox{resp.}\qquad \sum_{k=n+3}^{\infty} \,Y_{\a,k}$$
with the notations of the Lemma, so that one has  
$$\fa(x)\; =\; \int_0^\infty f_{\a,n}(x-y) g_{\a,n}(y) dy.$$
By the Proposition, the function $f_{\a,n}$ has a $\C^{n+1}$ extension on $\rl$ and one obtains 
$$\fa^{(n+1)}(x)\; =\; \int_0^\infty f_{\a,n}^{(n+1)}(x-y) g_{\a,n}(y) dy\; =\; \int_0^\infty f_{\a,n}^{(n+1)}(y) g_{\a,n}(x-y) dy.$$
The Laplace transform of $g_{\a,n}$ reads
$$\int_0^\infty e^{-\lb x} g_{\a,n}(x) \, dx\; =\; \prod_{k=n+3}^{\infty}\lpa \frac{1}{1+(\sin(\pi\a)/k\pi)^{1/\a}\lb}\rpa$$
and its reciprocal is of the P\'olya-Laguerre class $\E_1^*$ defined in \cite{K} p. 336. By Theorem 7.3.2. (b) in \cite{K} p. 345, this shows that the kernel $g_{\a,n}(x-y)$ is ${\rm TP}_\infty$ on $\rl\times\rl.$ Notice that $g_{\a,n}(x-y)$ is however not ${\rm STP}_\infty,$ because of the indicator function. Using the notation (3.1) p. 20 in \cite{K} with $I=(0,+\infty),$ Theorem 3.1. (a) p. 21 in \cite{K} and the Proposition entail
$$S^-(\fa^{(n+1)})\; \le\; S^-(f_{\a,n}^{(n+1)})\; =\; n+1.$$
Now since $\fa^{(n)}$ has isolated zeroes on $(0,+\infty)$ and vanishes at zero and infinity, Rolle's theorem yields
$$S^+(\fa^{(n)})\; \le\; S^-(\fa^{(n+1)})-1\; \le\; n$$
where we have used the notation (3.2) p. 21 in \cite{K} with $I=(0,+\infty).$ This shows that $\fa^{(n)}$ vanishes at most $n$ times, and hence exactly $n$ times, on $(0,+\infty).$ The proof is complete.

\fin

\section{Remarks and open questions}

\subsection{The case when $1/\a$ is an integer} In this situation the bell-shape follows easily from the independent factorization
\begin{equation}
\label{Will}
X_{1/n}\; =\; n^{-n} \Gamma_{1/n}^{-1}\,\times\,\cdots\,\times\,  \Gamma_{(n-1)/n}^{-1},
\end{equation}
which was pointed out in \cite{Wi} and also basically shown in \cite{K} pp. 121-122. We discard the explicit case $n=2$ and set $g_n$ resp. $h_n$ for the density of $n^{-n}\Gamma_{1/n}^{-1}$ resp. $\Gamma_{2/n}^{-1}\times\ldots\times\Gamma_{(n-1)/n}^{-1}.$ The multiplicative convolution yields
$$f_{1/n}(x)\; =\; \int_0^\infty g_n(xy^{-1})h_n(y)\frac{dy}{y}$$
and one can clearly differentiate under the integral. This entails
$$f_{1/n}^{(i)}(x)\; =\; \int_0^\infty g_n^{(i)}(xy^{-1})h_n(y)\frac{dy}{y^{i+1}}\; =\; \frac{1}{x^i}\int_0^\infty h_n(xy^{-1}) g_n^{(i)}(y) y^{i-1} dy$$
for every $x > 0, i\ge 1.$ The kernel $h_n(xy^{-1})$ is the composition of $(n-2)$ kernels that are ${\rm STP}_\infty$ on $(0, +\infty)\times(0,+\infty)$ by Theorem 2.1 p. 18 in \cite{K}, and is hence itself ${\rm STP}_\infty$ on $(0, +\infty)\times(0,+\infty)$ by the Binet-Cauchy formula. With the above notation for $S^-$ and $S^+$, Theorem 3.1 (b) p. 21 in \cite{K} and the aforementioned bell-shape property of $g_n$ yield
$$S^+(f_{1/n}^{(i)})\; \le\; S^-(g_n^{(i)})\; =\; i$$
for every $i\ge 1$ so that $f_{1/n}$ is bell-shaped, too. The above factorization (\ref{Will}) shows also easily that $\log X_{1/n}$ has a bell-shaped density on $\rl,$ a property which we believe to be true for all $\a\in (0,1).$ When $\a$ is rational, another multiplicative factorization of $\Xa$ in terms of inverse Beta and inverse Gamma random variables was obtained in Lemma 2 of \cite{S1}. But the latter seems useless for the bell-shape of $\Xa$ or $\log\Xa,$ because inverse or log Beta laws have neither bell-shaped densities nor ${\rm TP}_\infty$ convolution kernels. 

\subsection{Positive self-decomposable densities}  A positive random variable $X$ is self-decomposable when its Laplace transform reads
$$\EE[e^{-\lb X}]\; =\; \exp -\lcr\gamma_0\lb \,+ \, \int_0^\infty (1-e^{-\lb u}) \frac{k(u)}{u} du\rcr, \qquad\lb \ge 0$$
for some $\gamma_0\ge 0$ and $k : (0, +\infty)\to\rl^+$ non-increasing. For example, the positive stable law is self-decomposable with $\gamma_0 = 0$ and $k(x) = \a/\Gamma(1-\a) x^\a.$ 
It is known and easy to see - see Theorem 27.7 in \cite{S} - that $X$ has a density $f$ when $k\neq 0.$ The unimodality of $f$ was shown by Wolfe \cite{Wo}. In the following we will suppose that $\gamma_0 = 0$  w.l.o.g. and set $f_k$ for the positive self-decomposable density associated with the non-degenerate spectral function $k.$

When $k(0+) =+\infty,$ Theorem 28.4 (ii)  in \cite{S} and Lemma 2.5 in \cite{SY} show that $f_k$ is smooth and that all its derivatives vanish at zero and at infinity, and Theorem 1.4 in \cite{SY} entails that $f_k$ is strictly unimodal. When $k(0+) < +\infty,$ it is also known - see Remark 28.6 in \cite{S} - that $f_k$ is not smooth at zero and hence not bell-shaped. In view of the main result of the present paper, one might raise the

\begin{Conju} The density $f_k$ is bell-shaped if and only if $k(0+) =+\infty.$
\end{Conju}

Apart from the positive stable, examples of densities verifying the above conjecture are the inverse Gamma, or more generally the densities of $\Gamma_t^{-a}$ with $t > 0, a \ge 1.$ Indeed, the latter are self-decomposable by the HCM criterion \cite{B}, and bell-shaped by a direct computation. Other self-decomposable examples are infinite exponential sums, whose bell-shape property can be proved exactly in the same manner as for $\fa.$ As pointed out in the introduction, when $k(0+) =+\infty$ Theorem 5.1.(ii) in \cite{SY} already shows that $f_k^{(n)}$ vanishes only once before the first zero of $f_k^{(n-1)},$ for every $n\ge 2.$ The proof of this latter property relies on a repeated use of Steutel's equation - see e.g. Theorem 51.1 of \cite{S}. Analyzing the further zeroes of $f_k^{(n)}$ with the same method seems however not quite obvious, because of the memory involved in this integro-differential equation. 

Let us now consider the case where $k(0+)$ is finite. When $k(0+)<1,$ Theorem 1.3 in \cite{SY} shows that $f_k$ is decreasing. When $k(0+) =1,$ the situation is complicated because of the three different types ${\rm I}_2,{\rm I}_3$ and ${\rm I}_4$ of \cite{SY}. In particular $f_k$ may have a non-trivial modal interval - see Lemma 53.2 in \cite{S}. When $k(0+) >n+1$ for some $n\ge 0,$ then Theorem 28.4 (i)  in \cite{S} and Lemma 2.5 in \cite{SY} show that $f_k$ is $\C^n$ with derivatives vanishing at zero and at infinity, and again Theorem 1.4 in \cite{SY} entails that $f_k$ is strictly unimodal. The following is hence quite natural.

\begin{Conjd} If $k(0+) >n+1$ for some $n\ge 0,$ then $f_k^{(i)}\sim\{a_i\}$ for all $i = 0\ldots n.$
\end{Conjd}

An example supporting this conjecture is the density $f_{a,b}$ of $-\log \beta_{a,b},$ where $\beta_{a,b}$ is the Beta random variable with parameters $a>0, b >1.$ Indeed, one can show by a direct computation that $f_{a,b}$ is self-decomposable with spectral function
$$k_{a,b}(x)\; =\; \frac{e^{-ax}(1-e^{-bx})}{(1-e^{-x})}$$
so that $k_{a,b}(0+) = b,$ and the explicit formula 
$$f_{a,b} (x) \; =\; \frac{\Gamma(a+b)}{\Gamma(a)\Gamma(b)}e^{-ax}(1-e^{-x})^{b-1}$$
entails that $f_{a,b}$ satisfies the required properties.

\subsection{Two-sided stable densities} As mentioned in the introduction, stable densities are all visually bell-shaped on Holt and Crow's graphics, and it is an interesting problem to show the property rigorously. It is known that two-sided stable densities are real-analytic on $\rl$, never vanish, and that all their derivatives tend to zero at infinity - see Remarks 14.8 and 28.8 in \cite{S}. Hence, their $n$-th derivative vanishes at least $n$ times on $\rl$ by Rolle's theorem (a more complicated proof of this fact is also given in \cite{G1} pp. 234-235). The strict unimodality follows from Yamazato's theorem \cite{Y1} by analyticity. The article \cite{G1} had claimed the full bell-shape property for all stable densities with an argument relying on a certain integral representation through a Student-type kernel, and a multiplicative TP property. The integral representation is correct and basically amounts to Bochner's subordination and Zolotarev's duality - see Lemma 2 in \cite{G1}. However the multiplicative TP property claimed in Lemma 1 (iv) of \cite{G1} and crucial to bound the number of zeroes - see pp. 236-237 in \cite{G1}, is false in general. Indeed, Schoenberg's theorem shows - see pp. 121-122 and p. 390 in \cite{K} - that the kernel $K_\a(x,y)=\fa(e^{x-y})$ is ${\rm STP}_\infty$ if and only if $1/\a$ is the reciprocal of an integer. Actually, $K_\a$ is not even ${\rm TP}_2$ for $\a > 1/2$ - see the main result of \cite{S1}. In a forthcoming paper, we also show that $K_\a$ is not ${\rm TP}_n$ as soon as $\a > 1/n$ and $1/\a$ is not the reciprocal of an integer. In \cite{G1} the skewed Cauchy case is also treated separately via Zolotarev's representation - see (2.8) therein. All in all the main result of \cite{G1} must hence be reduced to the

\begin{Theo}[Gawronski] Two-sided $\a-$stable densities are bell-shaped when $\a =1,2$ or when $1/\a$ is an integer.
\end{Theo}

Notice that this result does not extend directly to the one-sided case in spite of the affirmation made in (iv) p. 239 of \cite{G1}, because the number of zeroes of a function sequence might clearly increase at the continuous limit. In order to cover all stable densities, the multiplicative method of \cite{G1} seems not appropriate because of the poor variation-diminishing properties of $K_\a$ - see however \cite{S2} for a somehow related multiplicative point of view, which works at the unimodal level. On the other hand, Yamazato's additive factorization extends to {\em all} two-sided stable densities and delivers yet another proof of their unimodality - see again \cite{Y2} pp. 600-601. Unfortunately the complete monotonicity of the factor $\Ya$ which is crucial in our argument to get the full bell-shape - see the above Remark 2 (a), is clearly lost in the two-sided case. \\

\noindent
{\bf Acknowledgements.} Ce travail a b\'en\'efici\'e d'une aide de l'Agence Nationale de la Recherche portant la r\'ef\'erence ANR-09-BLAN-0084-01.

\end{document}